\documentclass[12pt]{amsart}

\usepackage{amsmath} 
\usepackage{amsfonts}
\usepackage{amssymb}
\usepackage{mathrsfs}
\usepackage{asymptote}
\usepackage[top=1.25in, bottom=1in, left=1in, right=1in]{geometry}
\usepackage{etoolbox}
\usepackage{caption}
\usepackage[T1]{fontenc}
\newtheorem{theorem}[subsection]{Theorem}

\newtheorem{lemma}[subsection]{Lemma}

\makeatletter

\renewcommand\subsubsection{\@secnumfont}{\bfseries}%
\renewcommand\subsubsection{\@startsection{subsubsection}{3}
  \z@{.5\linespacing\@plus.7\linespacing}{-.5em}%
  {\normalfont\bfseries\itshape}}
  
\makeatother

\renewcommand{\epsilon}{\varepsilon}

\makeatletter
    \newcommand{\wc}{\@ifstar{\wcstar}{\wcnostar}}
\makeatother

\newcommand{\wcnostar}{\rightharpoonup}
\newcommand{\wcstar}{\overset{\ast}{\rightharpoonup}}

\usepackage{indentfirst}
\title{Solution to the Number Rotation Puzzle}
\author{Thomas Lam}
\address{Carnegie Mellon University tjlam@andrew.cmu.edu}
\date{}

\begin{document}

\maketitle

\begin{abstract}
\noindent The Number Rotation Puzzle (NRP) is a combination puzzle in which the goal is to rearrange a scrambled rectangular grid of numbers back into order via moves that consist of rotating square blocks of numbers of fixed size. Over all possible boards and rotating block sizes, we find all solvable initial configurations and provide algorithms to solve such configurations.  For sufficiently large board and rotating block sizes, solvability conditions depend only on parity restrictions, with special additional conditions for smaller sizes.  One special case leads to a novel construction of the exotic outer automorphism on $S_6$.  
\end{abstract}

\noindent\textit{2010 Mathematics Subject Classification}. Primary 05E99; Secondary 00A08

\noindent\textit{Key Words. combination puzzle, three-cycle, exotic automorphism}

\section{Introduction}

\begin{figure}
    \centering
    \captionsetup{justification=centering}
    
    \includegraphics[height=1in]{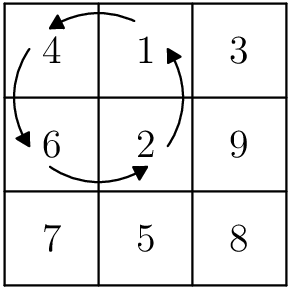}
    \includegraphics[height=1in]{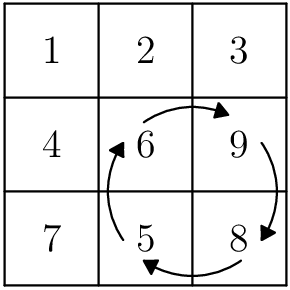}
    \includegraphics[height=1in]{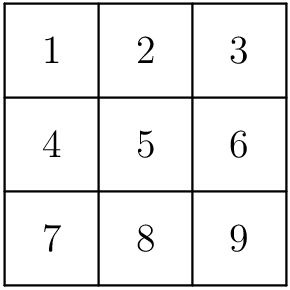}

    \caption{A two move solution to a sample standard NRP.}
    \label{fig:intro}
\end{figure}

The Number Rotation Puzzle (NRP) is a combination puzzle in which the goal is to rearrange a scrambled square grid of integers back into ascending order.  Numbers can be moved by rotating square blocks of numbers by some multiple of $90^\circ$.  The game is most commonly played on a $3 \times 3$ board of the integers $1$ to $9$, with $2 \times 2$ rotating blocks.  We call this the \textit{standard NRP}.  See Figure \ref{fig:intro} for an example of a game played on the standard NRP.

%%%  FIGURE 1  %%%
%\begin{wrapfigure}{R}{0.5\textwidth}
    %\fbox{
    %\begin{minipage}{3in}

    %\end{minipage}
    %}
%\end{wrapfigure}

In addition to the standard NRP, other variations exist, such as those with large board size and rotating block size, e.g. a $5 \times 5$ grid with $3 \times 3$ rotating blocks, and those with rectangular boards.

The NRP has been implemented on Nokia phones and a few puzzle game websites, most notably Simon Tatham's Twiddle (see \cite{twiddle}).
There has been very little research done on the NRP.  In fact, it has no ``official" name.  Most researchers focus on computational methods for solving the standard NRP.  For example, Lee programmed an AI using decision matrices to solve the standard NRP \cite{ai}, and Wang and Song developed an efficient brute-force algorithm to find solutions to the standard NRP \cite{bruteforce}.  However, one of the more outstanding feats in the research of the NRP is Fernando's solution to the $n \times n$ NRP with $(n-1) \times (n-1)$ rotating blocks, which was in turn used to bound God's Number --- the least number of moves required to solve a combination puzzle --- of the NRP \cite{ravi}.

The focus of this paper will be solving the NRP on rectangular boards, and thus generalizing the work of Fernando.  By solve, it is meant that all possible solvable initial configurations are identified, and a logical solving algorithm is developed that will solve all such configurations to prove that they are indeed solvable.  

Particularly, our main result is as follows:

\begin{theorem}
Consider the $n \times m$ board with $b \times b$ rotating blocks, where $n$ is the number of rows and $m$ is the number of columns.

    \begin{itemize}
    % \item[]
    \item If $b = n = m$, then there are exactly four solvable initial configurations, each a rotation of the solved board.  Assume in the following criteria that this case does not occur.
    \item If $b=1$, then the only solvable initial configuration is the solved board.
    \item If $b = 2$, then there are two further cases.
    \begin{itemize}
        \item If $(m,n) = (2,3),(3,2)$, then there are only 120 solvable configurations.
        \item Otherwise, all initial configurations are solvable.
    \end{itemize}
    \item If $b=3$, then there are two further cases.
    \begin{itemize}
        \item If $(m,n) = (3,4),(4,3)$, then there are only $6!$ solvable initial configurations.  In fact, solving one parity will solve the other.
        \item Otherwise, all initial configurations whose numbers lie in their correct parity set are solvable.
    \end{itemize}
    \item If $b \geq 4$, then case on the residue of $b$ modulo 8.
    \begin{itemize}
        \item If $b \equiv 2,6 \pmod{8}$, all initial configurations are solvable.
        \item If $b \equiv 0,4 \pmod{8}$, the solvable boards are those that are even permutations of the solved boards.
        \item If $b \equiv 3,\ 5 \pmod{8}$, the solvable boards are those that are even permutations of the solved boards such that each number lies in its correct parity.
        \item If $b \equiv 1,\ 7 \pmod{8}$, the solvable boards are such that each parity is an even permutation of their respective solved states.
    \end{itemize}
\end{itemize}
\end{theorem}

We will also reveal a surprising construction of the exotic outer automorphism on $S_6$ in the special case $(n,m,b) = (3,4,3)$.

Through this research, it is hoped that this puzzle be popularized in light of its intuitive mechanics.

\section{Preliminaries}

We denote by $(n,m,b)$ the NRP on an $n \times m$ board with $b \times b$ rotating blocks, where $n$ denotes the number of rows and $m$ denotes the number of columns.

We denote moves by capital letters, and typically such a move is a $90^\circ$ counter-clockwise rotation of a block.  If $X$ is such a move, then $X^{-1}$ is the $90^\circ$ clockwise rotation of said block, and $X^2$ is the $180^\circ$ rotation of said block.  We read algorithms from left to right, so that e.g. $XYZ$ denotes executing $X$, $Y$, then $Z$.  For an algorithm $A$ we write $A^n$ to denote $A$ repeated $n$ times, and write $A^{-1}$ to denote the inverse algorithm of $A$.

We also impose a coordinate system with $(i,j)$ denoting the square in the $i$th row from the top and the $j$th column from the left.  

\section{Parity Restrictions}

We let the case in which $m,n,b \geq 5$ with $m \neq n$ be the \textit{general case}.  We start by proving that our claimed conditions for solvability are necessary in the general case.  These conditions come from three \textit{parity restrictions} (PRs):

\begin{enumerate}
    \item Parity of a square
    \item Even permutations
    \item Even permutations within a parity of squares
\end{enumerate}

The ``parity of a square" PR simply refers to the invariance of the color of a square in which a number lies if the board is colored like a checkerboard.  Clearly this invariance occurs exactly when $b$ is odd, hence for odd $b$ we have that each of the numbers in the board must lie in the correct square parity in order to have solvability.  This gives us the first PR.

For some values of $b$, moves will execute even permutations, which indicates that solvable boards must be even permutations of the solved state.  To identify such $b$, we take cases:
\begin{itemize}
    \item If $b$ is even, then a move will move $b^2$ numbers in $b^2/4$ 4-cycles.  Since 4-cycles are odd, we have that the permutation introduced by a move is even iff $b^2/4$ is even, and this occurs exactly when $b \equiv 0 \pmod{4}$.
    \item If $b$ is odd, then each move will move $b^2-1$ numbers because the center of rotation stays in place.  Again, it follows that moves execute even permutations iff $\frac{b^2-1}{4}$ is even, and some simple casework shows that this occurs for all odd $b$. 
\end{itemize}
Hence for $b \equiv 0,1,3 \pmod{4}$, solvable boards must be even permutations, giving the second PR.

The third PR is a combination of the first two.  For odd $b$, each move executes a permutation that can be decomposed into a permutation of one parity and a permutation of the other parity.  If these permutations are even, then solvable boards must have each of their parities be even permutations of their respective solved states.  

Since a $b \times b$ rotating block moves $\frac{b^2-1}{2}$ numbers within each parity, we see as before that a move is an even permutation within a parity iff $\frac{b^2-1}{8}$ is even.  It is not hard to show that this occurs exactly when $b \equiv 1,7 \pmod{8}$, giving the third PR and completing the proof of the necessary direction in the general case.

\section{General Solving Algorithm}

We now show that these conditions are sufficient for solvability in the general case by constructing a solving algorithm.  We claim that it is sufficient to prove that
\begin{itemize}
    \item for even $n$, we may 3-cycle any three numbers, and
    \item for odd $n$, we may 3-cycle any three numbers lying in a common parity.
\end{itemize}

To see this, let us consider each of the cases described in the PRs.
\begin{itemize}
    \item If $n \equiv 0 \pmod{4}$, then the board must be an even permutation, and since we can execute any 3-cycle, we clearly can solve the board.
    \item For $n \equiv 2 \pmod{4}$, the board may be an odd permutation.  In this case, simply rotating any block $90^\circ$ will bring the board to an even permutation, in which case we can solve the board since we can execute any 3-cycle.
    \item For $n \equiv 1,7 \pmod{8}$, each parity must be an even permutation of their respective solved states.  Since we can compute any 3-cycle of numbers of the same number, it follows that we can execute any even permutation of a parity, hence the solvability.
    \item For $n \equiv 3,5 \pmod{8}$, the board must be an even permutation.  Looking at each parity, they are either both even permutations or both odd permutations of their respective solved states.  If the latter, then rotating any block $90^\circ$ must bring both parities to even permutations, from which we deduce solvability as in the previous case.
\end{itemize}

In practice, a solving algorithm would involve first executing a move if needed to bring the board to an even permutation, and then using 3-cycles to solve one number at a time until the entire board is inevitably solved.

Next, we claim that it is sufficient to solve the case in which $m = n+1$ and $b = n$, i.e. the $(n,n+1,n)$ NRP for $n \geq 5$.  One way to argue this is as follows:  For the $(M,N,n)$ NRP with $M 
\geq n$ and $N \geq n+1$, we may 3-cycle any three numbers for $n$ odd, and 3-cycle any three numbers on the same parity for $n$ even, by executing 3-cycles within the $n \times (n+1)$ sub-boards in which each of the three numbers lie until they lie in a common $n \times (n+1)$ sub-board.  Then we may 3-cycle them and undo the previous intermediate moves, resulting in the desired 3-cycle.

For the $(n,n+1,b)$ NRP, there are only two rotating blocks that may be used.  We denote rotation of the left and right rotating blocks counter-clockwise $90^\circ$ by $X$ and $Y$, respectively.  To prove that any three numbers may be 3-cycled, subject to parity restrictions, we will first prove that there exists an algorithm cycling some three numbers.

% Our rough sketch for the overall solving algorithm is as follows:
% \begin{enumerate}
%     \item Find a 3-cycle algorithm $\varphi$.
%     \item Pick three numbers to 3-cycle.  
%     \item Bring these numbers to the squares cycled by $\varphi$ in a series of intermediate moves.
%     \item Execute $\varphi$ and undo the intermediate moves.
%     \item This 3-cycles the desired numbers.  Repeat until solved.
% \end{enumerate}

\subsection{Existence of a 3-cycle Algorithm}

\begin{figure}
    \centering
    \includegraphics[height=2in]{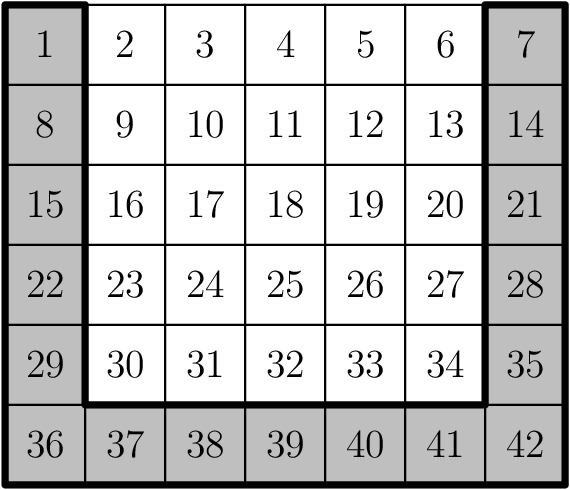}
    \includegraphics[height=2in]{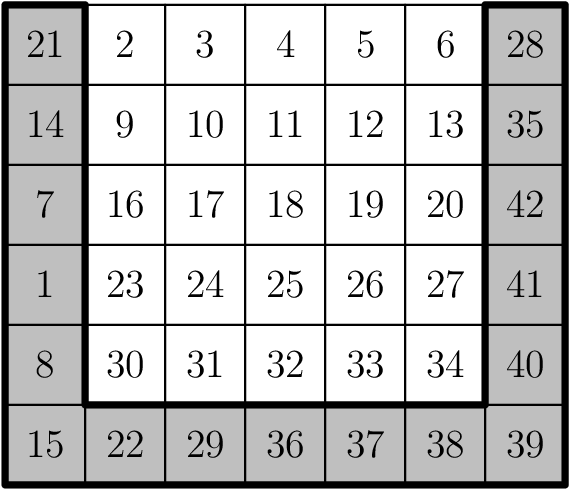}

    \caption{The algorithm \textbf{cycle} in action on the $6 \times 7$ board.  Left:  Initial board.  Right:  The board after execution of $(YX)^4$.  Since $n=6$, all numbers in the belt (outlined in bold) are cycled $6-3 = 3$ squares counter-clockwise.}
    \label{fig:cycle}
\end{figure}

We define an extremely useful algorithm, called \textbf{cycle} and denoted by $C$, as
$$C := (YX)^4 = YXYXYXYX.$$
The importance of $C$ is that it does not move many numbers.  Imagine that the numbers along the leftmost column, the bottom-most row, and the rightmost column form a sort of looped ``conveyor belt", with $(1,n+1)$ looping back to $(1,1)$.  Then $C$ will ``cycle" the numbers on the conveyor belt counter-clockwise a total of $n-3$ squares along the belt (See Figure \ref{fig:cycle}).  

To see that $C$ indeed executes this sort of behavior, first note that the numbers not in the belt form a square which is rotated $90^\circ$ with every execution of $YX$, hence these numbers' positions are unaffected after \textbf{cycle}.  Then, observe that the cyclic order of the numbers on the belt is invariant.  After a $YX$, the number at $(1,1)$ will move to $(n,1)$.  Thus, the belt is advanced $n-1$ squares counter-clockwise after each $YX$.  It follows that $(YX)^4$ advances the belt a total of $4n-4$ squares counter-clockwise.  Since there are only $3n-1$ squares on the belt, this is a net advancement of $n-3$ squares counter-clockwise.

We are now ready to define an algorithm $\varphi$ that executes a 3-cycle.  If we let $A = XYYXY^{-1}X^{-1}$, then we take
$$\varphi := C(ACA^{-1})C^{-1}(AC^{-1}A^{-1}).$$
The motivation is that the sequence $ACA^{-1}$ essentially ``applies \textbf{cycle} to a different set of squares", so in some sense $\varphi$ is simply commuting two different instances of \textbf{cycle}.  If the squares moved by $C$ and $ACA^{-1}$ intersect at exactly one square, then it is not hard to see that $\varphi$ must be a 3-cycle.

\begin{figure}
    \centering
    \includegraphics[height=2in]{cycle1.png}
    \includegraphics[height=2in]{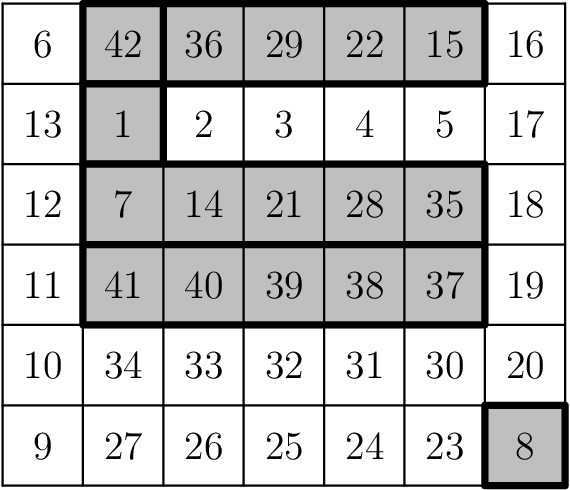}
    
    \caption{The cycle ``belt", after executing $A$.}
    \label{fig:A}

\end{figure}

To prove this, it is sufficient to show that after $A$, only one number that was in the belt will remain in the belt.  This can be visually seen in Figure \ref{fig:A}.  See Section A.1 for a complete proof.

% \textcolor{red}{\textbf{Is this a good compromise, or should I actually write out a full proof that this works?  }}

% A rigorous proof would be unsightly, so we opt to roughly explain the logic of $A$.  The motivation is to try and stack the first column, last column, and last row against each other and gather them up in the top rows.  See Figure \ref{fig:A} for the end result.

% \begin{itemize}
%     \item After $X$, the first column is sent to the bottom row, and the bottom row is sent to the second-to-last column, and is hence stacked up against the last column.
%     \item After $Y^2$, the image of the ``belt" is entirely contained in the first three columns and the top row.
%     \item $XY^{-1}$ moves these 
% \end{itemize}

% A rigorous proof would be unsightly.  See Figure \ref{fig:A}.  By inspection, the image of the ``belt" after executing $A$ is contained within the first rows minus the first and last columns, with the exception of the square in the bottom-right corner.  Since $n \geq 5$, this square is the sole element in the intersect of the ``belt" and its image after $A$, as needed.

\subsection{The Spiral Algorithm}

To build up our plan for sending the three numbers that we wish to 3-cycle to the squares cycled by $\varphi$, we first develop a methodology for sending one number to any desired square (while obeying PRs in the sense that the number's square and the goal square lie in the same parity when $n$ is odd).  To do this, it is sufficient to be able to send any number to a ``center square" of the board.  In the discussion that follows, we will assume that if $n$ is odd then the numbers we are maneuvering lie in the ``even parity" in the sense that their coordinates have even sum.  The case where $n$ is even and the numbers in question have odd parity is handled by reflecting our arguments.  

Our specific goal is to send a number, which we will call $a$, to the square $u=(n/2,n/2+1)$ if $n$ is even (the ``upper center square"), and $u=\left(\frac{n+1}{2},\frac{n+1}{2}\right)$ if $n$ is odd (the center of the $X$-rotating block).  Note that these squares are of minimum positive distance to the center of the $Y$-rotating block in either case.  The procedure for accomplishing this is the \textbf{spiral} algorithm, described as follows.

\begin{enumerate}
    \item If $a$ is located in the leftmost column, execute $X^2$ so that it is not.
    \item If either $a$ is at $u$, or is at any of the other three squares of minimum positive distance to the center of rotation of $Y$, then execute $Y$ as needed to move $a$ to $u$ and terminate.
    \item Execute $Y$ until the following inequalities hold, where $(i,j)$ are the coordinates of $a$:
    $$\frac{n+1}{2} \leq i \leq n, \qquad 2 \leq j \leq \frac{n+3}{2}.$$
    Visually, this condition entails that $a$ is located in the lower-left quadrant of the square rotated by $Y$, so this step is necessarily possible.
    \item Execute $X$.
    \item Loop back to Step 2.
\end{enumerate}

\begin{figure}
    \centering
    \includegraphics[height=2in]{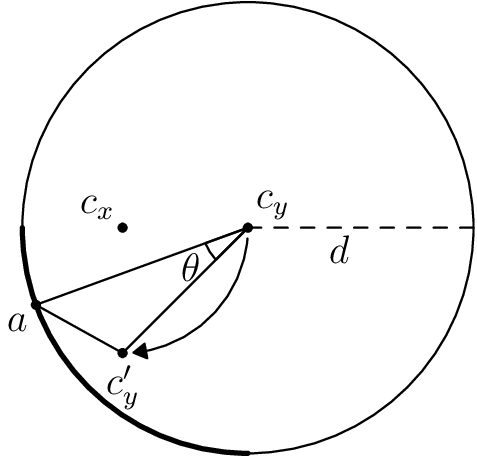}
    \caption{Geometrical interpretation of the \textbf{spiral} distance monovariant.  $c_x$ is the center of rotation of the move $X$, and similarly for $c_y$.}
    \label{fig:geo}
\end{figure}

We claim that this algorithm must terminate.  Let $c_y$ be the point at the center of rotation of $Y$, which may or may not be the center of a square.  It suffices to show that the straightline distance from $a$ to $c$ is a strictly decreasing monovariant over each loop of the algorithm, since there are only a finite number of possible distances between the centers of two squares.  Note that this distance is invariant as $Y$ is executed, so it suffices to observe a decrease in the distance after $X$ is executed.

Right before $X$ is executed, let the distance from $a$ to $c_y$ be $d$.  Then, since the coordinates of $a$ will satisfy the inequalities specified in step 3, $a$ will lie on the lower-left quarter arc of the circle centered at $c$ with radius $d$ (See Figure \ref{fig:geo}).  

Now, let us calculate a lower bound for $d$.  If $n$ odd, we see that the coordinates of $c_y$ sum to $n+1$, which is odd, so the premise that $a$ is at $c_y$ contradicts the assumption that $a$ lies on an even-parity square.  Moreover if $d=1$ then $a$ is already as close as it can be to $c_y$, which would have terminated the algorithm.  Hence we may assume $d > 1$ for odd $n$.  For even $n$, we see that the closest possible straightline distance to $c_y$ is $\frac{\sqrt{2}}{2}$, and if this distance is obtained then the algorithm would have terminated.  Otherwise, the next closest distance to $c_y$ would be $\frac{\sqrt{10}}{2}$, so we may assume that $d \geq \frac{\sqrt{10}}{2}$ for even $n$.  In either case, we may assume that $d > 1$.

Now, executing $X$ will rotate the aforementioned arc $90^\circ$ counter-clockwise about the center of rotation of $x$, denoted by $c_x$.  This will bring the arc closer to $c_y$.  Since only the distances between this arc and $c_y$ are of concern here, it is equivalent to rotate $c_y$ $90^\circ$ clockwise about $c_x$ to $c^\prime_y$.  The new distance between $a$ and $c_y$ after $X$ will simply be the distance between $a$ and $c_y^\prime$.  Let this new distance be $d^\prime$.

Now let the angle formed by $a$, $c_y$, $c^\prime_y$ be $\theta$.  Since $a$ must lie on the lower-left quarter arc, we see that $0^\circ \leq \theta \leq 45^\circ$.  The distance from $a$ to $c_y$ is the radius $d$, and since the distances between $c_x$ to $c_y$ and $c^\prime_y$ are both 1, the distance from $c_y$ to $c_y^\prime$ is $\sqrt{2}$.  By Law of Cosines we have that
$$d^\prime = \sqrt{d^2 + 2 - 2\sqrt{2}d\cos{\theta}}.$$
We want to show that $d^\prime < d$, or that
$$\sqrt{d^2 + 2 - 2\sqrt{2}d\cos{\theta}} < d.$$
After algebraic manipulations, we see that this holds if and only if $d\cos{\theta} > \frac{1}{\sqrt{2}}$.  This follows because we may multiply the assumption $d > 1$ with the inequality $\cos\theta \geq \frac{1}{\sqrt{2}}$ which holds for $0 \leq \theta \leq 45$.  This completes the proof that \textbf{spiral} terminates.  
%\textcolor{red}{\textbf{In retrospect, this is an arguably silly proof.  But, it's also pretty amusing!  Should I keep it or use a more ``sane" proof?}}  
We conclude that any number $a$ may be sent to any target square, provided that the square in which $a$ lies shares the same parity as the target square in the case that $n$ is odd.

We note that in \textbf{spiral}, there was nothing special about using the lower-left quadrant of the $Y$-rotating block and then executing $X$.  It would have been just as effective to instead move $a$ to the upper-left quadrant of the $Y$-rotating block and then executing $X^{-1}$.  Denote this variant as \textbf{spiral*}.

\subsection{The Spiral-Cycle Algorithm}

We are now ready to describe the procedure for sending any three numbers to any three desired squares.  As before, if $n$ is odd then we will require that the squares occupied by the three numbers, as well as the three goal squares, all lie on the same parity.  In this case we will assume without loss of generality that these six squares lie on the even parity, as before.

It suffices to choose three fixed squares $u_1,u_2,u_3$ and show that we can move any three numbers $a_1,a_2,a_3$ to these squares, since we can find the moves $B_1$ that send the three numbers we need to cycle to $u_1,u_2,u_3$, find the moves $B_2$ that send the goal squares to $u_1,u_2,u_3$, and then execute $B_1B_2^{-1}$.  

Our choice for $u_1,u_2,u_3$ is as follows:
\begin{itemize}
    \item If $n$ is even, we choose $u_1 = (1,1)$, $u_2 = (2,1)$, and $u_3 = (\frac{n}{2}, \frac{n}{2}+1)$.
    \item If $n$ is odd, we choose $u_1 = (1,1)$, $u_2 = (3,1)$, and $u_3 = (\frac{n+1}{2}, \frac{n+1}{2})$.  Note that these squares are chosen to lie on the even parity.
\end{itemize}

The procedure for sending $a_1,a_2,a_3$ to $u_1,u_2,u_3$ is the \textbf{spiral-cycle} algorithm, described as follows for all $n \geq 6$:

\begin{enumerate}
    \item Use \textbf{spiral} to move $a_1$ to $u_1$.
    \item If $a_2$ is now in the first column, we move it out via $XYX^{-1}$.
    \item Use a modified \textbf{spiral*} to bring $a_2$ to $(n,n+1)$ for even $n$ or $(n-1,n+1)$ for odd $n$, with the following caveats:  
    \begin{itemize}
        \item If we require usage of the $X$ rotating block in the \textbf{spiral*} algorithm, but rotating the $X$ block will displace the current position of $a_1$ (i.e. $a_1$ lies in the first column), then we execute the inverse \textbf{cycle} $C^{-1}$.  This ``hides" $a_1$ in the last column.  Note that $C^{-1}$ cannot displace $a_2$, which would have to lie in the upper-left quadrant of the $Y$ rotating block when \textbf{spiral*} requires the use of the $X$ block, so $a_2$ cannot lie in the belt at this time.
        \item If we require usage of the $Y$ rotating block in the \textbf{spiral*} algorithm, but rotating the $Y$ block will displace the current position of $a_1$ (i.e. $a_1$ lies in the last column), then we execute $C$.  This cannot displace $a_2$ either, since at this stage it must be the case that at least one $X^{-1}$ move of the \textbf{spiral*} algorithm was executed, so $a_2$ is too close to the center of the $Y$ rotating block to lie in the belt at this time.
        \item Once \textbf{spiral*} is terminated in this way, we execute $C$ is needed to restore $a_1$ to its position at $u_1 = (1,1)$.
    \end{itemize}
    \item Execute $XY^{-1}X^{-1}$.  This will bring $a_2$ to $(2,1)$ for even $n$, or $(3,1)$ for odd $n$.  Now $a_1,a_2$ are at $u_1,u_2$.
    \item If $a_3$ is now in the first column, then move it out using $XY^2XYXY^{-1}X^2Y^{-1}XY^{-1}X^{-1}$ for odd $n$, or $XY^2X^{-1}Y^{-1}XY^{-1}X^{-1}$ for even $n$.  See Section A.2 for a proof that this works for $n \geq 5$.
    \item Mimic Step 3 in order to get $a_3$ to $u_3$.  This will still work because $C^{-1}$ always moves the numbers on $(1,1),(2,1),$ and $(3,1)$ in the last column.
\end{enumerate}

If $n \geq 5$, replace each instance of $C$ with $C^2$.

To see that the application of \textbf{cycle} always ``hides" the numbers on $(1,1),(2,1),$ and $(3,1)$ as described, recall that $C^{-1}$ moves the numbers along the belt $n-3$ squares clockwise.  For $n \geq 6$ we see that $3 \leq n-3 \leq n$, so it must move the first three numbers on the belt into the last $n$, i.e. the last column.  For $n=5$, we use two applications of \textbf{cycle} instead so that the numbers move $2(5)-6 = 4$ squares clockwise along the belt.  Since $3 \leq 4 \leq 5$, this gives the same result of moving $(1,1),(2,1),$ and $(3,1)$ into the last column.  This justifies the correctness of \textbf{spiral-cycle}.  We conclude that
\begin{itemize}
    \item for $n$ even, we may send any three numbers to any three squares, and
    \item for $n$ odd, we may send any three numbers lying in a common parity to any three squares of that parity.
\end{itemize}

Consequently,
\begin{itemize}
    \item for $n$ even, we may 3-cycle any three numbers and
    \item for $n$ odd, we may 3-cycle any three numbers lying in a common parity,
\end{itemize}
which is what we wanted to prove.

% It is now simple to prove the sufficiency of the PRs for solvability in the general case.
% \begin{itemize}
%     \item If $n \equiv 0 \pmod{4}$, then the board must be an even permutation, and since we can execute any 3-cycle, we clearly can solve the board.
%     \item For $n \equiv 2 \pmod{4}$, the board may be an odd permutation.  In this case, simply executing either $X$ or $Y$ will bring the board to an even permutation, in which case we can solve the board since we can execute any 3-cycle.
%     \item For $n \equiv 1,7 \pmod{8}$, each parity must be an even permutation of their respective solved states.  Since we can compute any 3-cycle of numbers of the same number, it follows that we can execute any even permutation of a parity, hence the solvability.
%     \item For $n \equiv 3,5 \pmod{8}$, the board must be an even permutation.  Looking at each parity, they are either both even permutations or both odd permutations of their respective solved states.  If the latter, then executing of either $X$ or $Y$ must bring both to even permutations, from which we deduce solvability as in the previous case.
% \end{itemize}

This completes the proof of the theorem in the general case.

\section{Special Cases}

Our general solution fails for some small cases.  This is due to the fact that the general 3-cycle algorithms described do not have enough space to execute the permutation that contains the 3-cycle.  As a result, we have some special cases for smaller values of $m$, $n$, and $b$, that will introduce further solvability conditions.

\subsection{$(2,n,2)$}

$n$ is not large enough for any of the general algorithms to work.  In fact, for $n=3$, general 3-cycle algorithms do not exist.  It was proven by Jaap Scherphuis that out of all $6!$ theoretically reachable permutations in the $(2,3,2)$ variation, only $5!$ are achievable \cite{twogen}.

However, when $n \geq 4$, there exist algorithms that switch two numbers.  Take any $4 \times 2$ sub-board, and label the possible moves from left to right $X$, $Y$, and $Z$.  Then the following algorithm will switch two numbers:
$$XYZ^{-1}Y^2X^{-1}Z^{-1}YZ^2Y^{-1}$$

\subsection{$(m, n, 2)$ where $m,n \geq 3$}

Although this is classified as a special case, this is the most common version of the NRP, especially $(3,3,2)$, which is the standard NRP.  Thus, many different resolutions to this variation are well-known.

There are no restrictions on solvability.  Consider a $3 \times 3$ sub-board, and let $X$ rotate the lower-left $2 \times 2$ block, $Y$ rotate the upper-right block, and $Z$ rotate the lower-right block.  The following algorithm switches two numbers:
$$XY^{-1}X^{-1}YZ$$

From here, it is easy to prove that all initial configurations are solvable.

\subsection{$(3,4,3)$}

%\textcolor{red}{\textbf{This is a very important subsection.}}

We claim that out of all $12!$ theoretically reachable permutations, only $6!$ are achievable.  Moreover, we claim that solving one parity will solve the other.

\begin{figure}
    \centering
    \includegraphics[height=1in]{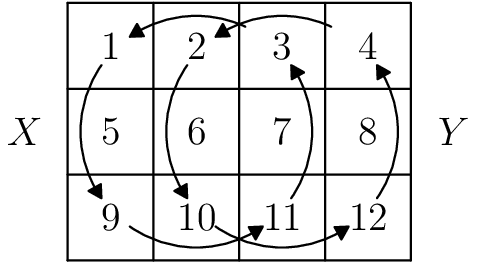}
    \caption{The $3 \times 4$ NRP with $3 \times 3$ rotating blocks.  The two possible moves are notated as shown.}
    \label{fig:5x3}
\end{figure}

To show this, first observe that since $b$ is odd, the parity in which each number lies is preserved by each move.  Hence we may view the puzzle as two separate but superimposed puzzles formed by by the 6 numbers in each parity. 

Now, let the two counter-clockwise moves be $X$ and $Y$, as in Figure \ref{fig:5x3}.  We let $P_0$ and $P_1$ be the set of coordinates with even coordinate sum and odd coordinate sum, respectively.  To prove our claim, we need to show that, if the initial configuration is solvable, then

\begin{enumerate}
    \item we can obtain any permutation of $P_0$ (without any regard as to what happens to $P_1$), and
    \item $P_0$ is solved if and only if $P_1$ is solved.
\end{enumerate}

To prove (1), note that $XYX^{-1}Y^{-1}XY^{-1}XYXYX$ swaps $(2,4)$ and $(3,3)$, and fixes all other elements of $P_0$.  It follows easily that any two numbers in $P_0$ may be swapped.  As for (2), it suffices to prove that $P_1$ being solved implies that $P_0$ is solved by symmetry.  In particular, it is sufficient to prove that any sequence of moves that fixes all elements of $P_1$ must also fix all elements of $P_0$.  This is because if for some solvable board we have that all elements of $P_1$ are solved but $P_0$ is not, then by virtue of solvability there exists an algorithm $A$ that solves the $P_0$ parity but fixes $P_1$.  What we wish to prove will contradict this.

We define a \textit{pair} to be a set of two elements.  Let $S_0$ be the set of all pairs with elements in $P_0$, and let $S_1$ be the set of partitions of $P_1$ into three pairs.  For a sequence of moves $A$, we let $A$ ``act" on $S_0$ and $S_1$ in an element-wise sense.  For example:
\begin{itemize}
    \item We have the pair $\{(1,1),(2,4)\} \in S_0$, and after executing the $Y^2$ this pair becomes $\{(1,1),(2,2)\}$.
    \item We have the partition $\{\{(1,2),(3,4)\},\{(1,4),(3,2)\},\{(2,1),(2,3)\}\} \in S_1$, and after executing $X$ this becomes the partition $\{\{(2,1),(3,4)\},\{(1,4),(2,3)\},\{(3,2),(1,2)\}\}$.
\end{itemize}
We shall denote such actions by $\cdot$, so e.g. $Y^2 \cdot \{(1,1),(2,4)\} = \{(1,1),(2,2)\}$.  %\textcolor{red}{\textbf{Is this notation terrible?}}

The remainder of the proof relies in the following remarkably strong property.

\begin{lemma}
There exists a bijection $\phi:S_0 \to S_1$ that preserves structure in the following sense:  For all pairs $p \in S_0$, we have that $\phi(A \cdot p) = A \cdot \phi(p)$ for all $A$.  
\end{lemma}
\label{lemma:graph}

\begin{figure}
    \centering
    \includegraphics[height=9.8cm]{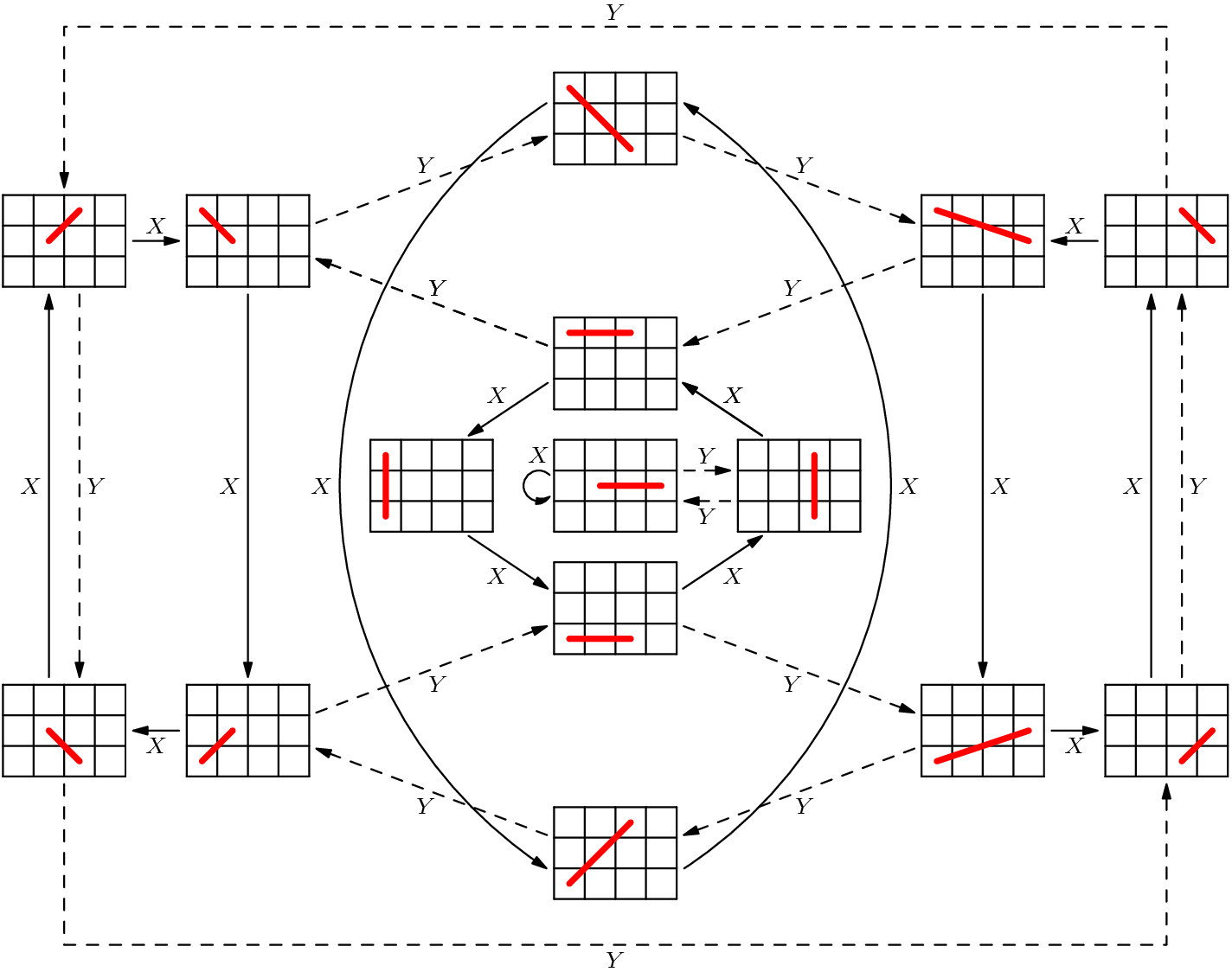}
    
    \includegraphics[height=9.8cm]{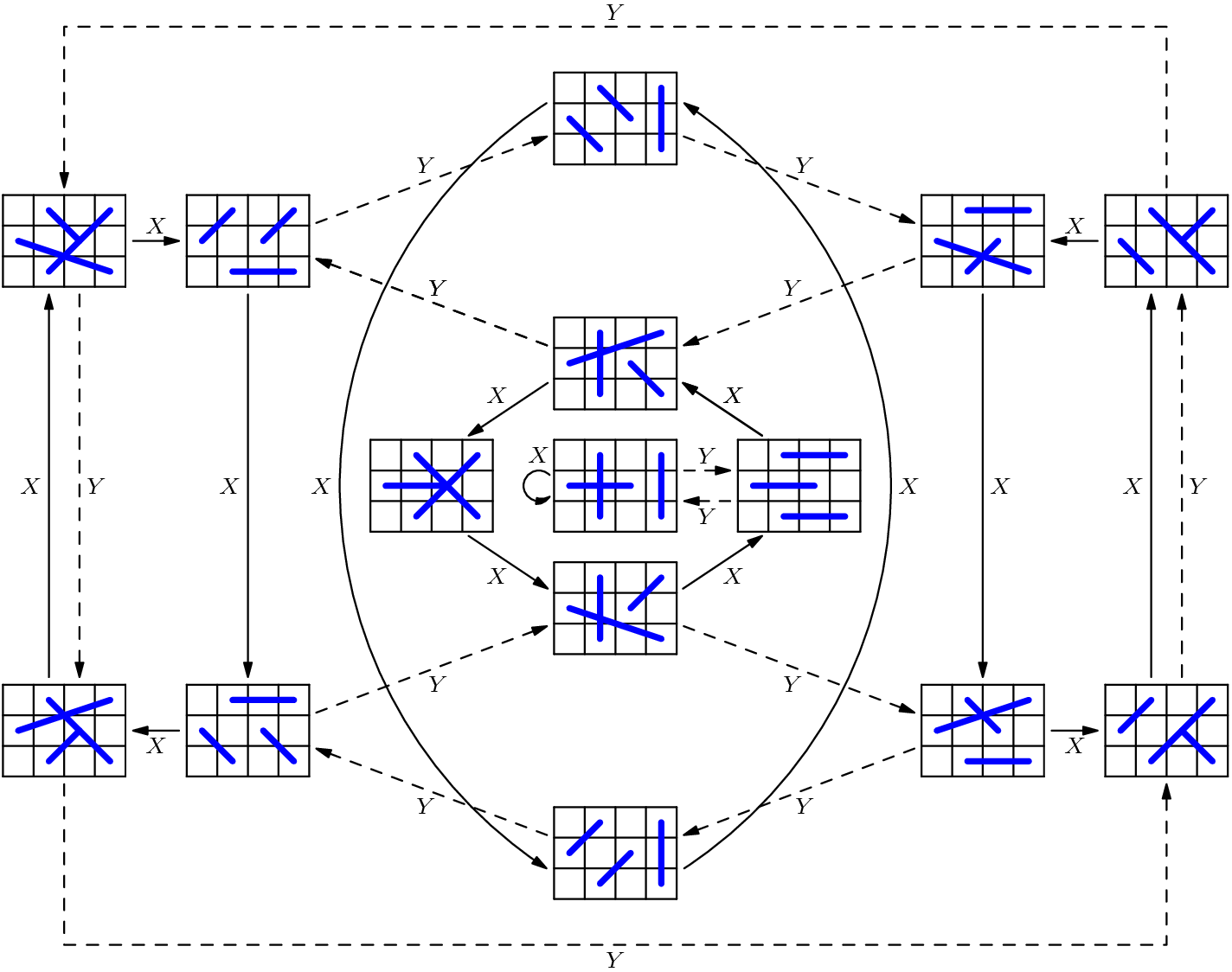}
    
    \caption{Movement graphs on the sets $S_1$ and $S_2$.  Top:  Movement graph on the set of pairs on the parity set $P_0$.  Bottom:  Movement graph on the set of partitions into pairs on the parity set $P_1$.  Each edge represents a move.  It can be seen that the structures of these graphs are identical.}
    \label{fig:graph}
\end{figure}

%Our proof of this is neither well-motivated nor necessarily elegant.  However, we will soon see a reason to believe that there is unlikely to be a nice proof.

\begin{proof}

% Observe that $|S_0| = \binom{6}{2} = 15$ and $|S_1| = \frac1{\displaystyle 3!}\binom{6}{2}\binom{4}{2}\binom{2}{2} = 15$, so $|S_0| = |S_1|$.   In fact, we claim that 

We list out all 15 elements of $S_0$ in a graph with two types of directed edges:  an ``$X$ edge" and a ``$Y$ edge", where each edge represents the actions of $X$ and $Y$.  That is, we place an $X$ edge between the pairs $p_1$ and $p_2$ if and only if $X \cdot p_1 = p_2$, and similarly for the $Y$ edges.  

We may construct the same sort of graph for $S_1$.  In doing so, we see that the two graphs are isomorphic!  See Figure \ref{fig:graph} for a visual of these graphs.  
%\textcolor{red}{\textbf{One improvement I plan to make is to rewrite this diagram in a program called Asymptote (so that it's more precise), and add arrows to the edges to indicate their direction.  I'll also label the arrows with an X or Y as described.}}  
The proof of our claim follows from this isomorphic structure by taking $\phi$ to be the natural correspondence that arises.

% We note the particular correspondence
% $$\phi(\{(1,1), (3,1)\}) = \{\{(2,1),(2,3)\},\{(1,2),(3,4)\},\{(3,2),(1,4)\}\},$$ 
% which can be seen from Figure \ref{fig:graph}.  From this, we claim that any sequence of moves that fixes the first column of squares $(1,1)$, $(2,1)$, and $(3,1)$ must also fix $(2,3)$, which is sufficient to prove the lemma.  Indeed, if a sequence of moves $S$ fixes $(1,1)$ and $(3,1)$, then $S \cdot \{(1,1), (3,1)\} = \{(1,1), (3,1)\}$, thus 
\end{proof}

To finish, consider a sequence of moves $A$ that fixes the elements of $P_1$.  Take any pair $\{a,b\} \in S_0$.  We claim that $A$ fixes this pair, i.e. $A \cdot \{a,b\} = \{a,b\}$.  Indeed, using the lemma and the fact that $A$ fixes the elements of $P_1$, we have that
$$A \cdot \{a,b\} = \phi^{-1}(A \cdot \phi(\{a,b\})) = \phi^{-1}(\phi(\{a,b\})) = \{a,b\}.$$
Since $A$ fixes all pairs, it must fix all elements of $P_0$.  This is because for any $a,b,c \in P_0$, we have $A \cdot \{a,b\} = \{a,b\}$ and $A \cdot \{a,c\} = \{a,c\}$, so $A \cdot a \in \{a,b\} \cap \{a,c\} = \{a\}$, forcing $A \cdot a = a$ for all $a \in P_0$, as needed.

As a remark, we note the following fascinating connection to group theory:  Let $\sigma_{P_0}$ and $\sigma_{P_1}$ denote that permutation groups on $P_0$ and $P_1$, respectively.
%\textcolor{red}{\textbf{Is this a standard notation?}}
Construct a map $\psi:\sigma_{P_0} \to \sigma_{P_1}$ as follows:  If a sequence of moves executes a permutation $\phi_0 \in \sigma_{P_0}$ on $P_0$ and a permutation $\phi_1 \in \sigma_{P_1}$ on $P_1$, then $\psi:\phi_0 \mapsto \phi_1$.  Then by the work we have done, $\psi$ is a well-defined isomorphism, and since $\sigma_{P_0} \cong \sigma_{p_1} \cong S_6$, we may view $\psi$ as an automorphism on $S_6$.  In fact, a quick check reveals that $\psi$ is the exotic outer automorphism on $S_6$!  

%\textcolor{red}{\textbf{1.  Should I add anything else?  2.  This neat exotic property was revealed to me in an email correspondence with Ravi Fernando.  Should I do anything to indicate this?}}

\subsection{$(3,n,3)$ for $n \geq 5$}
Consider some $3 \times 5$ sub-board.  Let $P_0$ and $P_1$ denote the squares in this sub-board of even and odd parity respectively, so that $P_0$ contains the upper-left corner.  Label, from left to right, the three possible moves, $X$, $Y$, and $Z$.  Then the algorithm $YZY^{-1}Z^{-1}YZ^2Y^{-1}Z^{-1}$ switches two numbers in $P_0$ (with other side effects in $P_1$), hence we can solve all numbers in $P_0$.  

Since each move is a rotation of a $3 \times 3$ block, each move enacts an even permutation, so the resulting configuration of the numbers in the other parity, $P_1$, must be an even permutation of the solved order.  Thus, it suffices to find an arbitrary 3-cycle algorithm on numbers in $P_1$ that does not displace the numbers in $P_0$.  Consider $(XZX^{-1}Z^{-1})^2$, which is a 3-cycle.  Then, one can demonstrate that any three numbers can be sent anywhere by adapting the transposition algorithm on $P_0$ to work on $P_1$ instead by "translating" it as such:  $XYX^{-1}Y^{-1}XY^2X^{-1}Y^{-1}$.  This will act as a transposition on $P_1$, and applying this repeatedly to move the three desired numbers to be cycled to the correct squares will serve as an intermediary algorithm.  

\subsection{$(m, n, 3)$ where $m, n \geq 4$}

The solvability condition here is the same as in the general case:  the initial configuration must be an even permutation of the solved board, and every number must lie in its correct parity.  Fernando's general 3-cycle algorithm on square boards is enough to resolve this case \cite{ravi}.  Nevertheless, we present a short 3-cycle algorithm.  Let $A$ rotate the upper-left $3 \times 3$ block, $B$ rotate the upper-right block, $C$ rotate the lower-left block, and $D$ rotate the lower-right block.  Then
$$ADA^{-1}D^{-1}C^{-1}DA^{-1}D^{-1}AC$$
is a 3-cycle.

\subsection{$(m,n,4)$ where $m,n \geq 4$ with $m \neq n$}

The solvability condition here is the same as in the general case:  the initial configuration must be an even permutation of the solved board.  As in the proof of the general case, it suffices to consider the case $(4,5,4)$ and find a 3-cycle.  An ``easy" 3-cycle is given by
$$(XY^2X^{-1}YX^{-1}Y^2)^{35}.$$

All remaining special cases are completely trivial.

\section{Future Research}

There are still some open problems concerning the NRP.  For example, all rotations were assumed to preserve the orientations of the numbers, so they always stay upright.  It would be interesting to determine the additional solvability conditions if orientation of numbers are not preserved with rotation.  Additionally, all moves were assumed to be rotations of square rotating blocks.  Thus, a variation of the NRP can be made by letting the rotating blocks be rectangular, with moves consisting of rotating said blocks $180^\circ$.  This would add an exciting dimension to an already convoluted puzzle.

\section{Acknowledgments}

I would like to thank Ravi Fernando for revealing to me the connection between my result on the $(3,4,3)$ NRP and the exotic automorphism on $S_6$.  

% \section{References}

\bibliographystyle{abbrv}
\bibliography{bibliography.bib}

\begin{thebibliography}{1}

\bibitem{ravi}
R.~Fernando.
\newblock Nxn corner rotation puzzle.
\newblock https://www.speedsolving.com/forum/threads/nxn-
  corner-rotation-puzzle.15472/, 2009.

\bibitem{ai}
Y.~Lee.
\newblock Solving the “rotation” puzzle in stages.
\newblock https://codemyroad.wordpress.com/2015/04/13
  /solving-the-rotation-puzzle-in-stages/, 2015.

\bibitem{twogen}
J.~Scherphuis.
\newblock Two-generator corners group.
\newblock URL:https://www.jaapsch.net/puzzles/pgl25.htm.

\bibitem{twiddle}
S.~Tatham.
\newblock Chapter 7: Twiddle, 2018.

\bibitem{bruteforce}
G.~Wang and J.~Song.
\newblock Bbfs-stt: An efficient algorithm for number rotation puzzle.
\newblock {\em Entertainment Computing}, 12:1--7, November 2015.

\end{thebibliography}

\newpage

\section{Appendix}

Here we include the unsightly proofs necessary to justify some of the claims we have made.

We recall the notations of $X$ and $Y$ denoting the $90^\circ$ counter-clockwise rotations of the left and right blocks respectively on the $(n,n+1,n)$ NRP.  We may view $X$ and $Y$ as functions on the board $\{1,\cdots,n\} \times \{1,\cdots,n+1\}$ defined via:
\begin{align*}
  X(x,y) &= \begin{cases}(n+1-y,x), & 1 \leq y \leq n \\ (x,y), & \text{otherwise}\end{cases}\\
  Y(x,y) &= \begin{cases}(n+2-y,x+1), & 2 \leq y \leq n+1 \\ (x,y), & \text{otherwise}\end{cases}  
\end{align*}
It will also be useful to state the behavior of $X^{-1},Y^{-1},X^2,$ and $Y^2$:
\begin{align*}
  X^{-1}(x,y) &= \begin{cases}(y,n+1-x), & 1 \leq y \leq n \\ (x,y), & \text{otherwise}\end{cases} \\
Y^{-1}(x,y) &= \begin{cases}(y-1,n+2-x), & 2 \leq y \leq n+1 \\ (x,y), & \text{otherwise}\end{cases} \\ 
X^2(x,y) &= \begin{cases}(n+1-x,n+1-y), & 1 \leq y \leq n \\ (x,y), & \text{otherwise}\end{cases}\\
Y^2(x,y) &= \begin{cases}(n+1-x,n+3-y), & 2 \leq y \leq n+1 \\ (x,y), & \text{otherwise}\end{cases}  
\end{align*}

For ease, let us define the abuse of notation $([a,b],c) := \{(a,c),(a+1,c),\cdots,(b,c)\}$, and likewise for $(a,[b,c])$.  For further ease and abuse, we view $\{(a,b)\}$ and $(a,b)$ as the same.

\subsection{Image of the belt under the algorithm $A$}

Let $A = XYYXY^{-1}X^{-1}$ be an algorithm, and let $E = ([1,n],1) \cup (n,[2,n]) \cup ([1,n],n+1)$ denote the ``belt".  We show that the image of $E$ under the algorithm $A$ intersects with $E$ at exactly one square, provided that $n \geq 5$.

We recall that our algorithms are read left-to-right, so eg. $(XY)(E) = Y(X(E))$.  Our proof strategy simply involves computing the image of $E$ after every move, and splitting said image to prepare for the next move.
\begin{itemize}
    \item After $X$:  
    \begin{align*}
    X(E) &= X([1,n],1) \cup X(n,[2,n]) \cup ([1,n],n+1) \\
     &= (n,[1,n]) \cup ([1,n-1],n) \cup ([1,n],n+1) \\
     &= (n,1) \cup (n,[2,n]) \cup ([1,n-1],n) \cup ([1,n],n+1)
    \end{align*}
    \item After $Y^2$:
    \begin{align*}
     (XY^2)(E) &= (n,1) \cup Y^2(n,[2,n]) \cup Y^2([1,n-1],n) \cup Y^2([1,n],n+1) \\
     &= (n,1) \cup (1,[3,n+1]) \cup ([2,n],3) \cup ([1,n],2) \\ 
     &= (n,1) \cup (1,n+1) \cup (1,[3,n]) \cup ([2,n],3) \cup ([1,n],2) 
    \end{align*}
    \item After $X$:
    \begin{align*}
        (XY^2X)(E) &= X(n,1) \cup (1,n+1) \cup X(1,[3,n]) \cup X([2,n],3) \cup X([1,n],2) \\
        &= (n,n) \cup (1,n+1) \cup ([1,n-2],1) \cup (n-2,[2,n]) \cup (n-1,[1,n])\\
        &= (n,n) \cup (1,n+1) \cup ([1,n-1],1) \cup (n-2,[2,n]) \cup (n-1,[2,n])\\
    \end{align*}
    \item After $Y^{-1}$:
    \begin{align*}
      &\phantom{{}={}} (XY^2XY^{-1})(E) \\
      &= Y^{-1}(n,n) \cup Y^{-1}(1,n+1) \cup ([1,n-1],1) \cup Y^{-1}(n-2,[2,n]) \cup Y^{-1}(n-1,[2,n])\\
      &= (n-1,2) \cup (n,n+1) \cup ([1,n-1],1) \cup ([1,n-1],4) \cup ([1,n-1],3)\\
    \end{align*}
    \item After $X^{-1}$:
    \begin{align*}
      &\phantom{{}={}} (XY^2XY^{-1}X^{-1})(E) \\
      &= X^{-1}(n-1,2) \cup X^{-1}([1,n-1],1) \cup X^{-1}([1,n-1],3) \cup X([1,n-1],4) \cup (n,n+1) \\
      &= (2,2) \cup (1,[2,n]) \cup (3,[2,n]) \cup (4,[2,n]) \cup (n,n+1) \\
    \end{align*}
\end{itemize}
From this we see that, under the key assumption that $n \geq 4$, we have that $(n,n+1)$ is the sole element in the intersection of $E$ and $A(E)$.

\subsection{Moving $a_3$ out of the first column}

After executing the first four steps of the Spiral-Cycle algorithm, it may be the case that $a_3$ lies in the first column, in which case we cannot proceed.  If so, then it must be moved out of this column without disturbing the positions of $a_1$ and $a_2$.  

% $XY^2XYXY^{-1}X^2Y^{-1}XY^{-1}X^{-1}$ for odd $n$, or $XY^2X^{-1}Y^{-1}XY^{-1}X^{-1}$ for even $n$.

Let us first consider the case in which $n$ is even.  We claim that after $XY^2X^{-1}Y^{-1}XY^{-1}X^{-1}$, all numbers in the first column are moved out of it, without displacing $a_1$ and $a_2$, which are are $(1,1)$ and $(2,1)$ respectively.

It is easy to verify that $a_1$ and $a_2$ are fixed under this algorithm.
\begin{align*}
  (1,1) &\stackrel{X}\mapsto (n,1) \stackrel{Y^2}\mapsto (n,1) \stackrel{X^{-1}}\mapsto (1,1)\\
  &\stackrel{Y^{-1}}\mapsto (1,1) \stackrel{X}\mapsto (n,1) \stackrel{Y^{-1}}\mapsto (n,1) \stackrel{X^{-1}}\mapsto (1,1) \\
  (2,1) &\stackrel{X}\mapsto (n,2) \stackrel{Y^2}\mapsto (1,n+1) \stackrel{X^{-1}}\mapsto (1,n+1)\\
  &\stackrel{Y^{-1}}\mapsto (n,n+1) \stackrel{X}\mapsto (n,n+1) \stackrel{Y^{-1}}\mapsto (n,2) \stackrel{X^{-1}}\mapsto (2,1)  
\end{align*}

Now let $E = ([3,n],1)$ be the rest of the first column.  We show that the image of $E$ under the algorithm is disjoint from $E$, which suffices.  Here we must assume $n \geq 4$.

\begin{itemize}
    \item After $X$:
    $$X(E) = (n,[3,n])$$
    \item After $Y^2$:
    $$(XY^2)(E) = Y^2(n,[3,n]) = (1,[3,n])$$
    \item After $X^{-1}$:
    $$(XY^2X^{-1})(E) = X^{-1}(1,[3,n]) = ([3,n],n)$$
    \item After $Y^{-1}$:
    $$(XY^2X^{-1}Y^{-1})(E) = Y^{-1}([3,n],n) = (n-1,[2,n-1])$$
    \item After $X$:
    $$(XY^2X^{-1}Y^{-1}X)(E) = X(n-1,[2,n-1]) = ([2,n-1],n-1)$$
    \item After $Y^{-1}$:
    $$(XY^2X^{-1}Y^{-1}XY^{-1})(E) = Y^{-1}([2,n-1],n-1) = (n-2,[3,n])$$
    \item After $X^{-1}$:
    $$(XY^2X^{-1}Y^{-1}XY^{-1}X^{-1})(E) = X^{-1}(n-2,[3,n]) = ([3,n],3)$$
\end{itemize}
Evidently $([3,n],3)$ is disjoint from $([3,n],1)$, as needed.

Now we consider the case in which $n$ is odd.  Then we claim that $XY^2XYXY^{-1}X^2Y^{-1}XY^{-1}X^{-1}$ moves out all numbers in the first column without displacing $a_1$ and $a_2$, which are at $(1,1)$ and $(3,1)$ respectively.

This algorithm is more difficult to motivate, but nevertheless we may manually verify that $a_1$ and $a_2$ are fixed under this algorithm.  
\begin{align*}
    (1,1) &\stackrel{X}\mapsto (n,1) \stackrel{Y^2}\mapsto (n,1) \stackrel{X}\mapsto (n,n) \stackrel{Y}\mapsto (2,n+1) \stackrel{X}\mapsto (2,n+1) \\
    &\stackrel{Y^{-1}}\mapsto (n,n) \stackrel{X^2}\mapsto (1,1) \stackrel{Y^{-1}}\mapsto (1,1) \stackrel{X}\mapsto (n,1) \stackrel{Y^{-1}}\mapsto (n,1) \stackrel{X^{-1}}\mapsto (1,1) \\
    (3,1) &\stackrel{X}\mapsto (n,3) \stackrel{Y^2}\mapsto (1,n) \stackrel{X}\mapsto (1,1) \stackrel{Y}\mapsto (1,1) \stackrel{X}\mapsto (n,1) \\
    &\stackrel{Y^{-1}}\mapsto (n,1) \stackrel{X^2}\mapsto (1,n) \stackrel{Y^{-1}}\mapsto (n-1,n+1) \stackrel{X}\mapsto (n-1,n+1) \stackrel{Y^{-1}}\mapsto (n,3) \stackrel{X^{-1}}\mapsto (3,1)
\end{align*}
The difficulty is showing that the rest of the column, $E = (2,1) \cup ([4,n],1)$, is moved out.  We will need the assumpton that $n \geq 5$.
\begin{itemize}
    \item After $X$:
    $$X(E) = X(2,1) \cup X([4,n],1) = (n,2) \cup (n,[4,n])$$
    \item After $Y^2$:
    $$(XY^2)(E) = Y^2(n,2) \cup Y^2(n,[4,n]) = (1,n+1) \cup (1,[3,n-1])$$
    \item After $X$:
    $$(XY^2X)(E) = X(1,n+1) \cup X(1,[3,n-1]) = (1,n+1) \cup ([2,n-2],1)$$
    \item After $Y$:
    $$(XY^2XY)(E) = Y(1,n+1) \cup Y([2,n-2],1) = (1,2) \cup ([2,n-2],1)$$
    \item After $X$:
    $$(XY^2XYX)(E) = X(1,2) \cup X([2,n-2],1) = (n-1,1) \cup (n,[2,n-2])$$
    \item After $Y^{-1}$:
    $$(XY^2XYXY^{-1})(E) = Y^{-1}(n-1,1) \cup Y^{-1}(n,[2,n-2]) = (n-1,1) \cup ([1,n-3],2)$$
    \item After $X^2$:
    \begin{align*}
        (XY^2XYXY^{-1}X^2)(E) &= X^2(n-1,1) \cup X^2([1,n-3],2) \\ 
        &= (2,n) \cup ([4,n],n-1)
    \end{align*}
    \item After $Y^{-1}$:
    \begin{align*}
        (XY^2XYXY^{-1}X^2Y^{-1})(E) &= Y^{-1}(2,n) \cup Y^{-1}([4,n],n-1) \\
        &= (n-1,n) \cup (n-2,[2,n-2])
    \end{align*}
    \item After $X$:
    \begin{align*}
        (XY^2XYXY^{-1}X^2Y^{-1}X)(E) &= X(n-1,n) \cup X(n-2,[2,n-2]) \\
        &= (1,n-1) \cup ([3,n-1],n-2)
    \end{align*}
    \item After $Y^{-1}$:
    \begin{align*}
        (XY^2XYXY^{-1}X^2Y^{-1}XY^{-1})(E) &= Y^{-1}(1,n-1) \cup Y^{-1}([3,n-1],n-2) \\
        &= (n-2,n+1) \cup (n-3,[3,n-1])
    \end{align*}
    \item After $X^{-1}$:
    \begin{align*}
        (XY^2XYXY^{-1}X^2Y^{-1}XY^{-1}X^{-1})(E) &= X^{-1}(n-2,n+1) \cup X^{-1}(n-3,[3,n-1]) \\ 
        &= (n-2,n+1) \cup ([3,n-1],4)
    \end{align*}
\end{itemize}
We see that $(n-2,n+1) \cup ([3,n-1],4)$ does not intersect the first column, as needed.

\end{document}